\documentclass{amsart}
\usepackage{amssymb}

\newtheorem{lemma}{Lemma}
\newtheorem{prop}[lemma]{Proposition}
\newtheorem{cor}[lemma]{Corollary}
\newtheorem{thm}[lemma]{Theorem}

\newtheorem{thm?}[lemma]{Theorem?}

\DeclareFontEncoding{OT2}{}{} 
  \newcommand{\textcyr}[1]{%
    {\fontencoding{OT2}\fontfamily{wncyr}\fontseries{m}\fontshape{n}%
     \selectfont #1}}
\newcommand{\Sha}{{\mbox{\textcyr{Sh}}}}

\title{The period-index problem in WC-groups I: elliptic curves}
\author{Pete L. Clark}
\address{1126 Burnside Hall \\ Department of Mathematics and Statistics \\
McGill University \\ 805 Sherbrooke West \\ Montreal, QC, Canada H3A 2K6}
\email{clark@math.mcgill.ca}

\newcommand{\F}{\ensuremath{\mathbb F}}
\newcommand{\Fp}{\ensuremath{\F_p}}

\newcommand{\Q}{\ensuremath{\mathbb Q}}

\newcommand{\Z}{\ensuremath{\mathbb Z}}

\newcommand{\ra}{\ensuremath{\rightarrow}}

\newcommand{\PP}{\ensuremath{\mathcal P}}

\newcommand{\Gal}{\operatorname{Gal}}
\newcommand{\Aut}{\operatorname{Aut}}

\newcommand{\Pic}{\operatorname{Pic}}
\newcommand{\FPic}{\underline{\Pic}}
\newcommand{\Gm}{\mathbb G_m}

\newcommand{\gk}{\mathfrak{g}_K}
\newcommand{\normresn}{\langle \ , \ \rangle_n} 
\newcommand{\normresp}{\langle \ , \ \rangle_p}
\newcommand{\normcc}{\langle C_1, C_2 \rangle_p}

\newcommand{\sep}{\operatorname{sep}}
\begin{document}
\maketitle

\begin{abstract}
Let $E/K$ be an elliptic curve defined over a number field, and let $p$
be a prime number such that $E(K)$ has full $p$-torsion.  We show that
the order of the $p$-part of the Shafarevich-Tate group of $E/L$ is unbounded
as $L$ varies over degree $p$ extensions of $K$.  The proof uses O'Neil's 
period-index obstruction.  We deduce the result
from the fact that, under the
same hypotheses, there exist infinitely many elements of the Weil-Ch\^atelet
group of $E/K$ of period $p$ and index $p^2$.
\end{abstract}

\noindent
\section{Introduction}
\noindent
The aim of this note is to prove the following result.
\begin{thm}
Let $p$ be a prime number, $E/K$ be an elliptic curve over a number field 
with $E[p](\overline{K}) = E[p](K)$, and $r$ be a
positive integer.  Then there are infinitely many degree $p$ field extensions
$L/K$ such that 
\[\dim_{\F_p} \Sha(E/L)[p] \geq r. \]
\end{thm}
\noindent
Recall that for any elliptic curve over a field $K$ of characteristic different from 
$p$, all $p$-torsion points become rational over an extension field of degree dividing
$|GL_2(\F_p)| = (p^2-1)(p^2-p)$.  Moreover, if $E/\overline{K}$ admits complex
multiplication,  all $p$-torsion points become rational over an extension of degree dividing
$2(p^2-1)$ or $2(p-1)^2$.  This immediately gives the following corollary.
\begin{cor}
If $E/K$ is an elliptic curve over a number field, $p$ a prime and $r$
a positive integer, there exist infinitely many field extensions $L/K$ of
degree at most $p^5$ such that $\dim_{\F_p} \Sha(E/L)[p] \geq r.$  Moreover,
for infinitely many $E/K$ -- namely those admitting complex multiplication
over $\overline{K}$ -- the same result holds for infinitely many field
extensions of degree at most $2p^3$.
\end{cor} 
\noindent
In Section 2 we deduce Theorem 1 as a consequence of the following result,
which is of independent interest.
\begin{thm}
Let $p$ be a prime, and $E/K$ an elliptic curve over a number
field with all its $p$-torsion defined.  Then there exists an infinite
subgroup of $H^1(K,E)[p]$ all of whose nonzero elements have index $p^2$
(i.e., $p^2$ divides the degree of any splitting field extension).
\end{thm}
\noindent
The proof of Theorem 3 makes essential use of the period-index obstruction map of Catherine 
O'Neil.
In particular we need two results concerning this map.  The first, Theorem 5,
gives a necessary and sufficient condition for the period to equal the index.
As the reader shall see, Theorem 5 is an immediate consequence of results of Cassels or 
of O'Neil.  The second, Theorem 6, is a computation
of the period-index obstruction in the case of full level structure, a result 
which appears in [O'Neil] but requires correction.  Section 3 is devoted to a review
of the period-index obstruction and a proof of these two results.
\\ \\
The proof of Theorem 3 is given in Section 4.
\\ \\
Finally, in Section 5 we discuss some issues raised by the proofs
and the possibility of certain generalizations.
\\ \\
Acknowledgements: It is a pleasure to acknowledge helpful conversations
with Catherine O'Neil and William Stein, in the course of which
many of the ideas of this paper emerged.
A motivation for the
particular form of the results presented here was the recent preprint
\cite{Kloosterman}, which proves -- by completely different methods -- a result 
similar to our Corollary 2 but involving
simultaneous variation of $E$ and $L/\Q$ among fields of degree $O(p^4)$.
Many thanks to the anonymous referee for several useful suggestions and one
invaluable remark.
\section{Theorem 3 implies Theorem 1}
\noindent
Let $S \subset H^1(K,E)[p]$ be
an infinite subgroup all of whose nontrivial elements 
have index $p^2$.  For each $\eta_i \in S$, there is a finite
set of places $v$ of $K$ such that $\eta_i$ remains nonzero
in the completion $K_v$.  By a theorem of [Lichtenbaum],
every class in $H^1(K_v,E)[p]$ can be split
by a degree $p$ extension.\footnote{For an alternate argument
using a theorem of Lang and Tate, see Remark 5.I.}  It now follows from
the usual (weak) approximation theorem for valuations that for any 
single $\eta_i$ we can find a degree $p$
global extension $K_i/K$ such that $\eta_i|_{K_i}$ is zero
everywhere locally, i.e., represents an element of 
$\Sha(E/K_i)[p]$.  Because the index of $\eta_i$ is $p^2$ and
we have made a field extension of degree only $p$, this is
a nontrivial element.  (This argument is due to William Stein.)
\\ \\
We now refine the above argument to produce $r$ $\Fp$-linearly
independent classes.  For this, observe first that $H^1(K_v,E)[p]$
is a finite group (e.g. it is a homomorphic image of $H^1(K_v,E[p])$,
and the Galois cohomology groups of a finite
module over a $p$-adic field are finite: \cite[Prop. II.5.14]{CG}).  
Starting with an element $\eta_1$ of
$S$, the subgroup $H_1 \subseteq S$ consisting of classes which are locally
trivial at all places where $\eta_1$ is locally nontrivial has finite
index and is therefore infinite; choose a nontrivial $\eta_2$ in this
group.  Continuing in this way, we can construct a cardinality
$r$ set $\{\eta_1,\ldots,\eta_r\}$ of $\Fp$-linearly independent elements of 
$S$ such that the sets $\Sigma_i$ of places where $\eta_i$ is locally
nontrivial are pairwise disjoint.  Accordingly, we can again find a single
global extension $L/K$ of degree $p$ such that all $r$ classes give
elements of $\Sha(E/L)[p]$.  Let $\eta = a_1 \eta_1 + \ldots + a_r \eta_r$
be any $\F_p$-linear combination of the $\eta_i$'s.  As above, if
$\eta|_L = 0$, then $\eta$ is a class in $S$ of index $p$, so
$\eta = 0$: i.e., $a_1 = \ldots = a_r = 0$.  
Thus $\dim_{\F_p} \Sha(E/L)[p] \geq r$.
\section{On the period-index obstruction for elliptic curves}
\noindent
Throughout this section the notation is
as follows: $K$ is an arbitrary field with
absolute Galois group $\mathfrak{g}_K = \Gal(K^{\sep}/K)$,
$n$ is a positive integer not divisible by the characteristic
of $K$, and $E/K$ is an elliptic curve.
\subsection{The period-index obstruction map}
Consider the \textbf{Kummer sequence}
\[ 0 \ra E(K)/nE(K) \ra H^1(K,E[n]) \ra H^1(K,E)[n] \ra 0. \]
The group $H^1(K,E)[n]$ parameterizes genus one curves
$C/K$ equipped with the structure of a principal homogeneous
space for $E = J(C) = \FPic^0(C)$ and having period dividing
$n$.  This geometric interpretation ``lifts'' to $H^1(K,E[n])$
as follows.
\begin{prop}
The group $H^1(K,E[n])$ classifies equivalence classes of
pairs $(C,[D])$, where $C$ is a principal homogeneous
space for $E$ and $[D] \in \FPic^n(C)(K)$ is a $K$-rational divisor \emph{class}
of degree $n$.  Two such classes are equivalent if and only
if there exists an isomorphism of principal homogeneous
spaces $f: C_1 \ra C_2$ such that $f^*([D_2]) = [D_1]$.
\end{prop}
\noindent
I have been unable to find this proposition in the literature
in the precise form in which we have stated it, but I am told that it has been well-known 
for a long time.  Indeed, Proposition 4 can readily be deduced either from work of Cassels or
of O'Neil.  
\\ \\
Sketch of proof: In either case, the idea is to interpret $E[n]$ as an automorphism group of 
a suitable structure
$\mathcal{S}$, so that by Galois descent $H^1(K,E[n])$ parameterizes the twisted
forms of $\mathcal{S}$.  But there is some latitude
in the choice of $\mathcal{S}$.  The classical choice 
\cite[Lemma 13.1]{Cassels2} is to view $E[n]$ as the deck 
transformation group of $[n]: E \ra E$, so that $H^1(K,E[n])$ parameterizes
finite \'etale maps $f: C \ra E$ which are geometrically Galois,
with group $\Z/n\Z \oplus \Z/n\Z$.
The correspondence is given as $f \mapsto (C,[nQ])$, where $Q$ is any element of
$f^{-1}(O)$.  O'Neil's choice [O'Neil, Prop. 2.2] is to view $E[n]$ as the automorphism 
group of the morphism $\varphi_{L(D)}: E \ra \PP^{n-1}$ associated to the ample divisor
$D = nO$, so that $H^1(K,E[n])$ parameterizes ``diagrams'' $C \ra V$, where $C$ is a principal
homogeneous space for $E$ and $V$ is a Severi-Brauer variety, indeed the Severi-Brauer variety
associated to an effective rational divisor class $D$ as in [CL, \S 10.6].  
\\ \\
Remark: O'Neil's method can be adapted to give an analogue of Proposition
4 for higher-dimensional abelian varieties; see \cite[\S 4.1]{Clark}.
\newcommand{\Spec}{\operatorname{Spec}}
\newcommand{\et}{\operatorname{\'et}}
\newcommand{\Br}{\operatorname{Br}}
\\ \\
Now recall that if $V/K$ is any 
smooth, projective, geometrically irreducible variety, there 
is an exact sequence
\[0 \ra \Pic(V) \ra \FPic(V)(K) \stackrel{\delta}{\ra} Br(K) \ra Br(V). \]
This is the exact sequence of terms of low degree arising from the Leray spectral sequence
\[E_2^{pq} = H^p(\Spec K, R^q f_* \Gm) \implies H^{p+q}(V,\Gm)\]
associated to the sheaf $(\Gm)_V$ and the morphism of \'etale sites induced
by $V \ra \Spec K$.
\\ \\
We define the 
\textbf{period-index obstruction map}
\[\Delta: H^1(K,E[n]) \ra Br(K) \]
by $(C,[D]) \mapsto \delta([D])$.  \\ \\
In terms of 
O'Neil's setup,
the obstruction map is simply given by 
\[(C \ra V) \mapsto [V] \in Br(K), \]
where by $[V]$ we mean the Brauer group element corresponding
to the Severi-Brauer variety $V$ [CL, \emph{loc. cit.}].
\begin{thm}
A class $\eta \in H^1(K,E)[n]$ of exact period $n$ has index $n$ if and only
if some lift of $\eta$  to $\xi \in H^1(K,E[n])$ has $\Delta(\xi) = 0$.
\end{thm}
\noindent
Proof: $\eta$ has index $n$ if and only if there exists a 
$K$-rational divisor of degree $n$ on the corresponding principal homogeneous
space $C$, i.e., if and only if some $K$-rational divisor class of
degree $n$ on $C$ has vanishing obstruction.  Thus the result is clear.
\subsection{Theta groups, Heisenberg groups and the explicit obstruction map}
Clearly Theorem 5 can only be useful if we know something about the
period-index obstruction map.  Happily, such knowledge is indeed
available, thanks to a result of O'Neil which identifies $\Delta: H^1(K,E[n]) 
\ra Br(K)$ as a connecting map in nonabelian Galois cohomology.  
\\ \\
Indeed, let $D = n[O]$ be the degree $n$ divisor on $E$ supported on the
origin, and let $L = L(D)$ be the associated line bundle.  Let
$\mathcal{G} = \mathcal{G}_L$ be the theta group associated to
the line bundle as in [Mumford]; $\mathcal{G}$ is an algebraic 
$K$-group scheme fitting into a short exact sequence
\begin{equation}
1 \ra \Gm \ra \mathcal{G} \ra E[n] \ra 1
\end{equation}
with center $Z(\mathcal{G}) = \Gm$.  Viewing (1) as a central
extension of $\gk$-modules, there is a 
cohomological connecting map 
\[\Delta: H^1(K,E[n]) \ra H^2(K,\Gm) = Br(K). \]
By [O'Neil, Prop. 2.3], $\Delta$ 
coincides with the period-index obstruction
defined in the last section.  
\\ \\
The goal of this section is to compute $\Delta$ in the case
when $E/K$ has full $n$-torsion defined over $K$.  That is,
we assume that
the finite \'etale $K$-group scheme $E[n]$ is \emph{constant},
and choose a Galois module isomorphism $E[n] \cong (\Z/n\Z)^2$. 
The Galois-equivariance of Weil's $e_n$-pairing implies that
$\Z/n\Z = \bigwedge^2 E[n] = \mu_n$ as Galois modules, so
the above choice of basis induces an isomorphism
\[H^1(K,E[n]) \cong H^1(K,\mu_n)^2 = (K^*/K^{*n})^2. \]
So in this case the period-index obstruction
can be viewed as a map
\[\Delta: (K^*/K^{*n})^2 \ra Br(K). \]
Now we must point out that [O'Neil, 3.4] gives a computation of $\Delta$
which is not quite correct: it is claimed that $\Delta(a,b) = \langle a, b 
\rangle_n$, the \textbf{norm-residue symbol}.  But the following counterexample was 
supplied by the referee:
\\ \\
Suppose $n=2$, so $E$ is given as $y^2 = (x-e_1)(x-e_2)(x-e_3)$.
Then the map $\iota: E(K)/2E(K)\ra (K^*/K^{*2})^2$ is given explicitly for
any point $(x,y) \in E(K)$ with $x \neq e_1, e_2$ as 
\[\iota(x,y) = (x-e_1,x-e_2) \pmod{K^{*2}} \]
[Silverman, Prop. X.1.4].  But $\Delta$
vanishes on $\iota(E(K)/2E(K))$, so in particular 
$\Delta(e_3-e_1,e_3-e_2) = 0$.  But as $e_1, \ e_2, \ e_3$
vary over all triples of distinct elements of $K$, $(e_3-e_1,e_3-e_2)$ runs
through all elements of $K^{\times}/K^{\times 2}$, and all Hilbert symbols
$\langle a, \ b \rangle_2$ vanish only if $Br(K)[2]$ vanishes.
\\ \\
On the other hand, the following result shows that the obstruction map
$\Delta$ is close to being the norm residue symbol $\normresn$.
\begin{thm}
Let $E/K$ be an elliptic curve over a field $K$ and $n$ a positive integer
not divisible by the characteristic of $K$ and such that $E[n]$ is
a trivial $\gk$-module.  Then there exist $C_1, \ C_2 \in K^*/K^{*n}$ such that
\[\Delta(a,b) = \langle C_1 a, C_2 b \rangle_n - \langle C_1, C_2 \rangle_n. \]
\end{thm}
\noindent
Before we begin the proof we will need to recall some facts about
Heisenberg groups.  There is an algebraic $K$-group scheme
$\mathcal{H}_n$, which is, like $\mathcal{G}$, a central extension of
$E[n]$ by $\Gm$.  To define $\mathcal{H}_n$, one chooses a decomposition
$E[n] = H_1 \oplus H_2$ into a direct sum of two cyclic order $n$ subgroup schemes.
(With a view towards the higher-dimensional case, one should think of this as a 
\textbf{Lagrangian decomposition}, i.e., that each $H_i$ is maximal isotropic for the 
Weil $e_n$-pairing; of course this is automatic for elliptic curves.)
Then $\mathcal{H}_n$ is defined by the following $2$-cocycle
$f_{H_1,H_2} \in Z^2(E[n],\Gm)$:
\[(P_1+P_2,Q_1+Q_2) \mapsto e_n(P_1,Q_2). \]
It is known (e.g. \cite[Prop. 2.3]{Sharifi}) that the coboundary map 
$\Delta_{\mathcal{H}}: H^1(K,E[n]) \ra Br(K)[n]$ associated to the Heisenberg
group $\mathcal{H}_n$ is nothing else but the norm-residue symbol $\langle \ , \
\rangle_n$.  Moreover, by a well-known result of [Mumford] the theta group
scheme $\mathcal{G}$ is isomorphic to the Heisenberg group $\mathcal{H}_n$
when the base field is separably closed.  Thus in general $\mathcal{G}$
is a Galois twisted form of $\mathcal{H}_n$.  Combining these two results 
with the above counterexample, it must be the case that
$\mathcal{G}$ can be a nontrivial twisted form of $\mathcal{H}_n$.
\\ \\
Nevertheless, we can completely understand
the possible twists: they are parameterized by
$H^1(K,\Aut_{\star}(\mathcal{H}_n))$, where the $\star$
indicates
that we want not the full automorphism group of
$\mathcal{H}_n$ but only the automorphisms which act trivially on the subgroup $\Gm$
and on the quotient $E[n]$.  It will turn out that $\Aut_{\star}(\mathcal{H}_n) 
\cong (H_1 \oplus H_2)$, so that the twisted forms of the Heisenberg group will
be parameterized by pairs of order $n$ characters of $\gk$.
\\ \\
We now begin the proof of Theorem 6.  Let
$\psi \in \Aut_{\star}(\mathcal{H}_n)$, and let 
$(P_1,P_2,\epsilon)$ denote
an arbitrary element of the Heisenberg group.  Since $\psi$ is the
identity modulo the center, we have $\psi(P_i) = P_i$ for $i = 1, \ 2$;
together with the fact that $\psi(0,0,\epsilon) = (0,0,\epsilon)$, this
implies that $\psi: (P_2,P_2,\epsilon) \mapsto 
(P_1,P_2,\chi(\psi)(P_1,P_2)\epsilon)$.  That is, an automorphism of
$\mathcal{H}_n$ as an extension determines a map $\chi: H_1 \oplus H_2 \ra \Gm$,
i.e., a \emph{character} of $H_1 \oplus H_2$.  Conversely, any such character
defines an automorphism, and we have canonically $\Aut_{\star}(\mathcal{H}_n) = 
(H_1 \oplus H_2)^{\vee}$ (Pontrjagin = Cartier dual).  It follows that
the collection of twisted forms of the Heisenberg group is $
H^1(K,(H_1 \oplus H_2)^{\vee}) \cong H^1(K,H_1 \oplus H_2)$, since
the Weil pairing gives an autoduality $E[n]^{\vee} \cong E[n]$.
\\ \\
Changing notation slightly, let \[\chi \in H^1(K,\Aut_{\star}(\mathcal{H}_n)) =
H^1(K,(H_1 \oplus H_2)^{\vee}) \cong (K^*/K^{*n})^2 \] be a one-cocycle.
Using $\chi$ we build a twisted form $\mathcal{H}_{\chi}$ of $\mathcal{H}_n$, i.e., 
the 
group scheme whose $\overline{K}$-points are the same
as the $\overline{K}$-points of $\mathcal{H}_n$, but with twisted $\gk$-action, as 
follows: \[\sigma \cdot (P_1,P_2,\epsilon) = (P_1,P_2,\chi(\sigma)(P_1,P_2)
\sigma(\epsilon)). \]
We may now compute the cohomological coboundary map $\Delta$ directly from its 
definition.  For this, we view
$\mathcal{H}_{\chi}/\overline{K}$ as $\Gm \times E[n]$ ``doubly twisted,''
i.e., twisted as a $\gk$-set as just discussed, and twisted as a group
via the cocycle $f$ introduced above:
\[(\alpha,P) \star (\beta,Q) = (\alpha \beta f(P,Q),P+Q). \]
We note that the inverse of $(\alpha,P)$ is $(\alpha^{-1}f(P,-P)^{-1},-P)$.
Let $\eta \in Z^1(K,E[n])$; we want to compute $\Delta(\eta)(\sigma,\tau)$.
The basic recipe for this allows us to choose
arbitrary lifts $N_{\sigma}, \ N_{\tau}, \ N_{\sigma \tau}$ of
$\eta_{\sigma}, \ \eta(\tau) \ \eta(\sigma \tau)$ to $\mathcal{H}_{\chi}$
and put $\Delta(\eta)(\sigma,\tau) = N_{\sigma}\sigma(N_{\tau})N_{\sigma \tau}^{-1}$.
We choose to lift by the set-theoretic identity section: $\eta(\sigma) \mapsto
(1,\eta(\sigma))$, and so on.  Keeping in mind that $\sigma(\eta(\tau)) = \eta(\tau)$
and $\eta(\sigma \tau) = \eta(\sigma) \eta(\tau)$, we get:

\[\Delta(\eta)(\sigma,\tau) = (1,\eta(\sigma)) \star \sigma(1,\eta(\tau)) \star
(1,\eta(\sigma \tau))^{-1} = \] 
\[(1,\eta(\sigma)) \star(\chi(\sigma)(\eta(\tau)),\eta(\tau)) \star 
(f(\eta(\sigma \tau),-\eta(\sigma \tau))^{-1},-\eta(\sigma \tau) = \]
\[(\chi(\sigma)(\eta(\tau))f(\eta(\sigma),\eta(\tau)),\eta(\sigma) \eta(\tau)) \star
(f(\eta(\sigma \tau),-\eta(\sigma \tau))^{-1},-\eta(\sigma \tau)) = \]
\[(\chi(\sigma)(\eta(\tau))f(\eta(\sigma),\eta(\tau)),0). \]
That is, the coboundary map $\Delta: H^1(K,E[n]) \ra Br(K)[n]$ is a product of two terms:
\[\Delta(\eta)(\sigma,\tau) = \Delta_1 \cdot \Delta_2 = \chi(\sigma)(\eta(\tau)) \cdot 
f(\eta(\sigma),
\eta(\tau)). \]
Indeed $\Delta_2$ and $\Delta_1$ are respectively the \emph{quadratic form} and
the \emph{linear form} comprising the quadratic map $\Delta$.  Both of these terms
are now easily recognizable: the quadratic part $\Delta_2$ (which, notice, is equal
to $\Delta$ when $\mathcal{H}_{\chi} = \mathcal{H}_n$) is the norm residue
symbol [Sharifi, \emph{loc. cit.}].
\\ \\
To evaluate $\Delta_1$, choose a basis $(P_1,P_2)$ of $E[n]$ and use the induced 
decomposition of $E[n] = H_1 \oplus H_2$
and the corresponding decomposition of the dual space
$E[n]^{\vee}$ (i.e., we decompose any character $\phi$ into $\psi_1 \oplus \psi_2$,
where $\chi_i(H_{j}) = 0$ for $i \neq j$).  This induces decompositions $\eta = 
\eta_1 \oplus \eta_2$ and $\chi = \chi_1 \oplus \chi_2$, so that
\[\chi(\sigma)(\eta(\tau)) = \chi_1(\sigma)(\eta_1(\tau))  \cdot \chi_2(\sigma)(\eta_2(\tau)).
 \]
Now under our identification $H^1(K,E[n]) = (K^{*}/K^{*n})^2$, $\eta_1$ corresponds
to $a \pmod{K^{* n}}$ and $\eta_{2}$ corresponds to $b \pmod{K^{* n}}$, so $\Delta_1$
is just the sum of the cyclic algebras $(a,\chi_1)$ and $(b,\chi_2)$.  Using
Kummer theory to identify the characters with elements (say) $C_2, \ C_1'$ of $K^*/K^{*n}$, we
get 
\[\Delta_1(a,b) = \langle a, C_2 \rangle + \langle b, C_1' \rangle =
\langle a, C_2 \rangle + \langle C_1, b \rangle, \]
where $C_1 = C_1'^{-1}$.
Thus we have 
\[\Delta(a,b) = \langle a, b \rangle + \langle a,C_2 \rangle + \langle C_1, b \rangle =
\langle C_1a,C_2b \rangle - \langle C_1,C_2 \rangle, \]
completing the proof of the theorem.
\section{The Proof of Theorem 3}
\noindent
In this section the following hypotheses are in force: $n = p$ is prime,
$K$ is a number field, and $E/K$ is an elliptic curve with 
$E[p](\overline{K}) = E[p](K)$.  We note that this implies,
by the Galois-equivariance of the Weil pairing, that
$K$ contains the $p$th roots of unity.  
Since for any class $\eta \in H^1(K,E)[p]$ the 
possible lifts of $\eta$ to $H^1(K,E[p])$ are parameterized by the
\emph{finite} abelian group $E(K)/pE(K)$ (weak Mordell-Weil theorem), by Theorem 5
the proof of Theorem 3 is reduced to the following result.
\begin{prop}
Let $K$ be a number field containing the
$p$th roots of unity
and $H \subseteq (K^{*}/K^{* p})^2$ a finite subgroup.  Then
there exists an infinite subgroup $G \subseteq (K^*/K^{* p})^2$
with the property that for every nonzero element $g$ of $G$ and  every element $h \in H$, 
$\Delta(hg) \neq 0$.
\end{prop}
\newcommand{\kp}{K^*/K^{*p}}
\noindent
By Theorem 6, $\Delta = \normresp$ up to a linear term, and essentially what must
be shown is the same statement with $\normresp$ in place of $\Delta$; this says,
morally, that Brauer groups of number fields are ``large'' in a certain sense.
We prove this directly (if inelegantly) using exactly what the reader expects:
local and global class field theory, especially the nondegeneracy of the local
norm residue symbol.
\\ \\
Along these lines we will need the following routine result, whose proof we include 
for 
completeness.
\begin{lemma}
Let $n$ be a positive integer, $K$ be a number field containing the
$n$th roots of unity, and $L_1,\ldots,L_k$ be $k$ cyclic degree $n$ extensions
of $K$.  Then the image in $K^*/K^{*n}$ of the subgroup of $K^*$ consisting
of simultaneous norms from each $L_i$ is infinite.
\end{lemma}
\noindent
Proof: By Hasse's norm theorem, if $L/K$ is a cyclic extension of
number fields, then $a \in K^*$ is a norm from $L$ if and only if it is
everywhere a local norm.  Let $S$ be the set of places of $K$ consisting of
the real Archimedean places (if any) together with all finite places which
ramify in any $L_i/K$ (if any).  Let $G_1 \subseteq K^*$ be the subgroup of elements which
are $n$th powers locally at every $v \in S$; notice that $G_1$ has finite index.
Recalling that the norm map on an unramified local 
extension is surjective onto the unit group, we get that any $a \in G_1$ is a simultaneous
local norm except possibly at the unramified places $v$ at which it has nontrivial valuation.
Let $h$ be the class number of $K$.  Then the set of primes which split completely
in the Hilbert class field as well as in each $L_i$ has density at least
$\frac{1}{hn^k}$.   For such a $v$, let $\pi_v$ be a generator of the corresponding prime ideal, 
and
let $G_2$ be the (infinite) subgroup of $K^*$ generated by these elements $\pi_v$.  Since
$G_1$ has finite index, $G := G_1 \cap G_2$ remains infinite and visibly has infinite
image in $K^*/K^{* n}$; by Hasse, every element of $G$ is a simultaneous norm.
\\ \\
Now we begin the proof of Proposition 7.  Write out the elements of $H$ as follows:
\[H = \{(h_{1i},h_{2i})\} | \ 1 \leq i \leq k\}.\]
Moreover, let $B = B_H$ be the finite set of places of $K$ containing the Archimedean places, the places
at which any $h_{1i}$ or $h_{2i}$ has nonzero valuation, and the places for which,
for any $e \pmod p$, any local symbol $e\langle C_1, C_2 \rangle_v - \langle h_{1i},h_{2i} \rangle_v$ is nonzero in $Br(K_v)$.
\\ \\
Clearly it's enough to construct arbitrarily large finite
subgroups $G$ such that every nontrivial element $(g_1,g_2)$ of $G$
has the property that for all $i$,
\[
\Delta(h_ig) = \langle C_1h_{1i}g_1, C_2h_{2i}g_2 \rangle_p \neq \normcc. \]

\noindent
We make two preliminary simplifying assumptions: first, let $C$ be the cyclic
subgroup generated by $\normcc$ in $Br(K)[p]$.  Rather than constructing elements
$g$ such that all modifications of $g$ by elements of $H$ have $\Delta(hg) \neq
\normcc$, it is convenient for a later inductive argument to require the stronger
property that for all $h \in H$, $\Delta(hg)$ is not an element of $C$.  Second,
by replacing $H$ by $H+C$, we reduce to the following problem: find
arbitrarily large finite subgroups $G$ all of whose nontrivial elements
$(g_1,g_2)$ have the property that for all $h = (h_{1i},h_{2i})$ in $H$,
\newcommand{\isnotin}{\ \operatorname{is} \ \operatorname{not} \ \operatorname{in} \ }
\begin{equation}
\langle h_{1i}g_1, h_{2i}g_2 \rangle_p \isnotin C.
\end{equation}
In order to accomplish this, we first claim that we can choose $g_2 \in \kp$ such that: \\ \\
$\bullet$ For $1 \leq i \leq k$, $\langle h_{1i},g_2\rangle  = 0$; and \\
$\bullet$ For $1 \leq i \leq k$, $g_2h_{2i} $ is not in $K^{*p}$. \\ \\
Indeed, the elements $g_2$ satisfying the first condition are precisely the
simultaneous norms from the $k$ cyclic field extensions  $K(h_{1i}^{1/p})/K$, so in the notation
of Lemma 6 there is a positive density set $S_1$ of principal prime ideals $v = (\pi_v)$ such 
that
$\pi_v \in K^*$ is a simultaneous norm from these $k$ extensions.  The second condition is
also satisfied as long as $v \in S_1 \setminus B$, so choose any such $v$ and take $g_2 = \pi_v$.
\\ \\
If we now choose any $g_1$ with the property that for all $i$ and any $e \pmod p$
\[\langle g_1,g_2h_{2i} \rangle  \neq e\normcc - \langle h_{1i}, h_{2i} \rangle, \]
then the element $g = (g_1,g_2)$ will have the desired property (2).  For each $i$,
since $g_2h_{2i}$ is not a $p$th power, there exists an infinite set of places $v = v(i)$
such that $g_2h_{2i}$ is not a $p$th power in $K_{v}$.  Hence we may choose
places $v_1, \ldots, v_k$, distinct and disjoint from $B$, such that for all $i$,
$g_2h_{2i}$ is not a $p$th power in $K_{v_i}$.  
By weak approximation, we can choose an element $g_1$ of
$\kp$ such that for all $i$, $g$ completes to a class of
$K_{v_i}^*/K_{v_i}^{*p}$ making all the local norm residue symbols $\langle g_1,g_2h_{2i} \rangle _{v_i}$
nontrivial (this is possible because of the nondegeneracy of the local
norm residue symbol).  But by definition of $B$, $e\langle C_1, C_2 \rangle_{v_i} - 
\langle h_{1i},h_{2i} \rangle_{v_i} = 0$ for all $i$, so 
we have constructed an element $g = (g_1,g_2)$ satisfying (2).
Now observe that if $1 \leq j < p$, $g_2^j$ satisfies
the same two bulleted properties as $g_2$; moreover, since $H$ is a subgroup,
$h_{2i}  = h_{2i'}^j$ for some other index $i'$, and the nontriviality of
$\langle g_1,g_2h_{2i} \rangle _v$ implies the nontriviality of $\langle g_1^j,g_2^j h_{2i}^j\rangle _v$,
so that indeed the entire cyclic subgroup $A$ generated by $(g_1,g_2) $ has the desired
property (2).
\\ \\
We finish by iterating the construction: running through the above
argument with $H$ replaced by $A \oplus C$ gives a two-dimensional
$\F_p$-subspace of $\kp \times \kp$, and so on.
\section{Concluding remarks}
\noindent
I. In the derivation of Theorem 1 from Theorem 3, instead of appealing to Lichtenbaum's
theorem on the equality of the period and index for \emph{all} classes in the Weil-Ch\^atelet 
group of an elliptic curve over a local field, we could instead have used
an earlier result of [Lang-Tate] giving the same equality for abelian varieties
of arbitrary dimension over local fields in the case when $p$ is prime to the residue
characteristic and $A$ has good reduction.  Indeed the set of places of $K$ lying over
$p$ together with those places of bad reduction for $E/K$ form a finite set, and as
in the proof we need only restrict to the finite index subgroup of classes trivial
at all these places.  
\\ \\
II. The proof of the main theorem shows that each nonzero element $g$ of $G \subseteq H^1(K,E)[p]$ gives
rise to at least one set of ``local conditions'' on a degree $p$ extension $L/K$ sufficient
to ensure that $g$ restricts to a nonzero element of $\Sha(E/L)$.  On the other hand, the proof
of Theorem 3 shows that $G$ is not only an infinite subgroup but has (in some sense)
``positive measure,'' bounded away from zero in terms of $\# E(K)/pE(K)$.  Thus the argument
should lead to an explicit lower bound on the function
\[f(N) = f(E/K,p,N) := \sum_{L/K,\ [L:K] = p,\ ||\Delta_{L/K}|| \leq N} \dim_{\F_p} \Sha(E/L)[p], \]
 where $\Delta_{L/K}$ is the discriminant of $L/K$.  What is to be expected about the
asymptotics of $f$?  
\\ \\
III. The hypothesis that $E$ has full $p$-torsion defined over $K$ is used only in the appeal
to the ``explicit'' period-index obstruction of Theorem 5 in the proof of Theorem 3.
My hope is that Theorem 3 should be valid for every elliptic curve over a number field
-- namely, there should always exist an infinite subgroup of principal homogeneous
spaces of order $p$ and index $p^2$.  The challenge
here is to make (sufficiently) explicit the period-index obstruction map $\Delta: H^1(K,E[p]) \ra
Br(K)$ in the case of an arbitrary Galois module structure on $E[p]$.  Notice that
the setup of Theorem 4 can be generalized to the case of elliptic curves $E$ such
that $E[n]$ has a Lagrangian decomposition: i.e., a decomposition
into one-dimensional subspaces $H_1 \oplus H_2$ as Galois module.  To be sure,
this is still quite a stringent condition -- satisfied, for any fixed non-CM
elliptic curve over a number field, for at most finitely many primes $p$
-- but at least this condition can be
satisfied for elliptic curves over $\Q$ for the primes $2$, $3$ and
$5$: for such primes, elliptic curves $E/\Q$ with Galois module structure 
$E[n] \cong \mu_p \oplus \Z/p\Z$ are known to exist.  In these cases, an analogue
of Theorem 4 would show the existence of genus $1$ curves $C/\Q$ of period
$p$ and index $p^2$ for $p \leq 5$.  While this may not sound very impressive,
we must point out that heretofore the \emph{only} examples in the literature
of genus one curves over \emph{any} number field with index exceeding their period are those of period $2$ 
and index $4$ (over $\Q$) constructed by \cite{Cassels} more than 40 years ago.  Indeed, Cassels' 
Jacobian elliptic curves have full $2$-torsion over $\Q$, so his results are a special
case of our Theorem 3.
\\ \\
IV. A question: for which fields $K$ is Proposition 7 valid?  Two obvious necessary
conditions are that $Br(K)$ be nontrivial and that the group of $p$th power classes 
$\kp$ be infinite.  Surely the proposition will hold
for all finitely generated fields of characteristic zero.  (I was so repulsed by the ugliness
of the proof that I have not even attempted such
a generalization.  I hope that someone else can find a cleaner way to proceed.)
Since the weak Mordell-Weil theorem 
remains valid
in this context \cite{Lang-Tate}, we would get a generalization of Theorem 3.
We note that the problem of finding necessary and sufficient
conditions on a field $K$ for there to exist genus one curves of period
$p$ and index $p^2$ remains open.  
\\ \\
IV$'$. Certainly there are \emph{some} examples of period-index violations in Weil-Ch\^atelet
groups of elliptic curves over function fields over number fields: let $K$ be a number field, 
$E_0/K$ an elliptic curve and $\eta_0 \in H^1(K,E_0)[p]$ a class of index $p^2$.  Let
$L := K(T_1,\ldots,T_n)$.  Write $E := E_0 \times_K L$ for the corresponding
(constant) elliptic curve over $L$; similarly put $\eta := \eta_0 |_{G_L}$.
Using the facts that $Br(K) \hookrightarrow Br(L)$ and $E(L) = Maps_K(\mathbb{P}^n,E) = 
E(K)$, we see that since $\eta_0$ does not admit a lift with
non-vanishing obstruction, neither does $\eta$.  But much more should be true:
if $L = K(V)$ is any function field of a variety over a number field and $E/L$
is any elliptic curve, then for all primes $p$ we expect $H^1(L,E)[p]$ to contain
infinitely many classes of index $p$ and infinitely many of index $p^2$.  
\\ \\
V. There are versions of Theorem 1 and Theorem 3 for principal homogeneous spaces
over abelian varieties of any dimension.  The proofs require a higher dimensional
analogue of the period-index obstruction and are pursued in a forthcoming paper
\cite{Clark}.

\end{document}